\documentclass[pdflatex,sn-mathphys-num]{sn-jnl}

\usepackage{amsmath,amssymb,amsfonts,amsthm,mathtools}
\usepackage{lmodern}
\usepackage{graphicx}
\usepackage{aliascnt}
\usepackage{enumitem}
\usepackage{cleveref}
\setlist[enumerate]{label=(\alph*),ref=(\alph*)}

\newcommand{\R}{\mathbb{R}}
\newcommand{\C}{\mathbb{C}}
\newcommand{\Lp}[1]{\mathrm{L}^{#1}}
\newcommand{\Hk}[1]{\mathrm{H}^{#1}}
\newcommand{\conC}{\mathrm{C}}

\newcommand{\X}{\mathcal{X}}
\newcommand{\Y}{\mathcal{Y}}
\newcommand{\U}{\mathcal{U}}
\newcommand{\V}{\mathcal{V}}
\newcommand{\W}{\mathcal{W}}
\newcommand{\Z}{\mathcal{Z}}
\newcommand{\Rf}{\mathbf{R}}

\newcommand{\Lb}{\mathcal{L}_{\mathrm b}}
\let\Re\relax

\DeclareMathOperator{\Re}{Re}

\DeclareMathOperator{\diag}{diag}
\DeclareMathOperator{\dom}{dom}
\DeclareMathOperator{\ran}{ran}

\DeclareMathOperator{\graph}{graph}

\newcommand{\adjunsymb}{\ast}
\newcommand{\adjunX}[1]{^{\adjunsymb_{#1}}}
\newcounter{adjunlevel}
\newtoks\adjunstars
\newcommand{\adjun}[1][1]{%
  \setcounter{adjunlevel}{1}%
  \adjunstars={\adjunsymb}%
  \loop\ifnum\value{adjunlevel}<#1%
    \adjunstars=\expandafter{\the\expandafter\adjunstars\adjunsymb}%
    \stepcounter{adjunlevel}%
  \repeat%
  ^{\the\adjunstars}%
}

\newcommand{\llangle}{\langle\!\langle}
\newcommand{\rrangle}{\rangle\!\rangle}
\DeclarePairedDelimiter{\norm}{\lVert}{\rVert}

\DeclarePairedDelimiterX{\setdef}[2]{\{}{\}}{#1\,\delimsize\vert\,\mathopen{} #2}
\DeclarePairedDelimiterX{\scprod}[2]{\langle}{\rangle}{#1,#2}
\DeclarePairedDelimiterX{\dualprod}[2]{\langle}{\rangle}{#1,#2}
\DeclarePairedDelimiterX{\sdprod}[2]{\llangle}{\rrangle}{#1,#2}
\newcommand{\pvek}[2]{\begin{pmatrix}#1\\#2\end{pmatrix}}
\newcommand{\spvek}[2]{\left(\begin{smallmatrix}#1\\#2\end{smallmatrix}\right)}

\newcommand{\sbmat}[4]{\left[\begin{smallmatrix}#1&#2\\#3&#4\end{smallmatrix}\right]}
\newcommand{\ddt}{\frac{\text{\normalfont d}}{\text{\normalfont d}t}}
\newcommand{\ddts}{\tfrac{\text{\normalfont d}}{\text{\normalfont d}t}}
\newcommand{\dd}{\mathrm{d}}
\newcommand{\dx}[1][x]{\mathop{\dd#1}}

\newcommand{\Flow}{\mathcal{F}}
\newcommand{\Effort}{\mathcal{E}}
\newcommand{\D}{\mathcal{D}}
\newcommand{\Lagr}{\mathcal{L}}
\newcommand{\Rc}{\mathcal{R}}

\newcommand{\Hsym}{H_{\mathrm{sym}}}

\theoremstyle{thmstyleone}
\newtheorem{theorem}{Theorem}[section]
\newaliascnt{proposition}{theorem}
\newtheorem{proposition}[proposition]{Proposition}
\aliascntresetthe{proposition}
\newaliascnt{lemma}{theorem}

\aliascntresetthe{lemma}
\newaliascnt{corollary}{theorem}
\newtheorem{corollary}[corollary]{Corollary}
\aliascntresetthe{corollary}
\theoremstyle{thmstyletwo}
\newaliascnt{example}{theorem}

\aliascntresetthe{example}
\newaliascnt{remark}{theorem}
\newtheorem{remark}[remark]{Remark}
\aliascntresetthe{remark}
\newcommand{\smalldiamond}{\raisebox{0.15ex}{\scalebox{0.8}{$\lozenge$}}}
\makeatletter
\newcommand{\inlineqed}[1]{%
  \leavevmode\unskip\penalty9999\hskip0pt\hbox{}\nobreak\hfill\hbox{#1}%
}
\let\oldremark\remark
\let\endoldremark\endremark
\renewenvironment{remark}{\oldremark}{\inlineqed{\smalldiamond}\endoldremark}
\makeatother
\theoremstyle{thmstylethree}
\newaliascnt{definition}{theorem}
\newtheorem{definition}[definition]{Definition}
\aliascntresetthe{definition}

\begin{document}
\title[Descriptor representations of geometric port-Hamiltonian systems]
{Descriptor representations of geometric port-Hamiltonian systems on Hilbert spaces}
\author[1]{\fnm{Anna} \sur{Fischer}}
\email{anna.fischer@tu-ilmenau.de}
\author*[1]{\fnm{Timo} \sur{Reis}}
\email{timo.reis@tu-ilmenau.de}
\affil*[1]{\orgdiv{Institut f\"ur Mathematik},
\orgname{Technische Universit\"at Ilmenau},
\orgaddress{\street{Weimarer Stra\ss e 25}, \postcode{98693}
\city{Ilmenau}, \country{Germany}}}

\abstract{We establish conversions in both directions between geometric and descriptor representations of linear port-Hamiltonian systems on infinite-dimensional Hilbert spaces. For the geometric-to-descriptor conversion, we represent the possibly multivalued inverse Dirac relation by a skew-adjoint operator on enlarged Hilbert spaces, retaining the original operator-domain restrictions. The construction accommodates internal and external ports. We show that the conversions preserve trajectories and the power balance. An incompressible Kelvin--Voigt solid illustrates the results.}
\maketitle

\section{Introduction}\label{sec:introduction}
Port-Hamiltonian systems separate energy storage, power-conserving
interconnection, and dissipation. This makes power exchange in coupled
systems transparent, particularly for distributed-parameter models
with boundary interactions
\cite{vanderSchaftMaschke2002,LeGorrecZwartMaschke2005,
JacobZwart2012,BrugnoliHaineMatignon2023}.
System-node formulations and boundary-triplet methods offer
operator-theoretic approaches to these infinite-dimensional systems
\cite{PhilippReisSchaller2025,GernandtPhilippPreusterSchaller2026}.

Linear port-Hamiltonian systems have two common formulations.
The \emph{geometric formulation} describes admissible trajectories
\emph{behaviorally}: at each $t$ in an interval $I$, they satisfy
the Dirac, Lagrangian, and resistive inclusions
\begin{equation}\label{eq:geometric-system}
 \left(
 \begin{pmatrix}
  \dot x(t)\\
  -f_{\mathrm R}(t)\\
  -y(t)
 \end{pmatrix},
 \begin{pmatrix}
  e_{\mathrm S}(t)\\
  e_{\mathrm R}(t)\\
  u(t)
 \end{pmatrix}
 \right)\in\D,\qquad
 (x(t),e_{\mathrm S}(t))\in\Lagr,\qquad
 (f_{\mathrm R}(t),e_{\mathrm R}(t))\in\Rc.
\end{equation}
Here $x$ is the state, its derivative $\dot x$ is the storage flow,
$(f_{\mathrm R},e_{\mathrm R})$ denotes the internal resistive port,
and $u$ and $y$ are the external input and output. The relations
specify which trajectories are admissible without first solving for
the state derivative or the output.

In the \emph{descriptor formulation}, the same dynamics are expressed
in an operator-based input--output form,
\begin{equation}\label{eq:compact-descriptor}
 \begin{pmatrix}
  \ddts\bigl( E\zeta(t)\bigr)\\
  -y(t)
 \end{pmatrix}
 = M 
 \begin{pmatrix}
  Q\zeta(t)\\
  u(t)
 \end{pmatrix}.
\end{equation}
The variable $\zeta$ collects storage, resistive, and auxiliary
variables; $E$ and $Q$ are bounded operators, while
$M$ is closed and may be unbounded. Its domain is part of the
equation and can impose additional restrictions. In this form the
external input $u$ and output $y$ are singled out, although existence
or uniqueness of a trajectory is not presupposed.

To relate this formulation to the established finite-dimensional
theory, consider the time-invariant port-Hamiltonian descriptor
systems of \cite[Def.~4]{BeattieMehrmannXuZwart2018},
\[
 \begin{aligned}
 \ddts\bigl(Ez(t)\bigr)&=(J-R)Qz(t)+(B-P)u(t),\\
 y(t)&=(B+P)\adjun Qz(t)+(S+N)u(t).
 \end{aligned}
\]
These equations can be written in the form
\eqref{eq:compact-descriptor} by setting
\[
 M=
 \begin{bmatrix}
  J-R&B-P\\
  -(B+P)\adjun&-(S+N)
 \end{bmatrix}.
\]
Here $J=-J\adjun$, $N=-N\adjun$, and the dissipation
matrix
$\left[\begin{smallmatrix}
  R&P\\
  P\adjun&S
 \end{smallmatrix}\right]$
is positive semidefinite. The energy matrices satisfy
$E\adjun Q=Q\adjun E\geq 0$.
In finite dimensions, $M$ can be split into a
skew-Hermitian and a negative semidefinite Hermitian part.
For unbounded operators $A$, this procedure may fail:
$\dom A\cap\dom(A\adjun)$ can be trivial, even if $A$
is maximal dissipative \cite{ArlinskiiTretter2020}.
We therefore derive the power structure directly from the
geometric relations. In infinite dimensions, $M$ may be
unbounded or non-densely defined; its domain is essential
and, together with auxiliary variables, can encode constraints.

For infinite-dimensional descriptor equations, see also
the operator-theoretic treatment in
\cite{ErbayJacobMorris2024}.
Our aim is to convert the geometric and descriptor
formulations into one another without changing the
admissible trajectories or their power balance.

The relation between geometric and descriptor formulations has
been studied extensively in finite dimensions
\cite{vanderSchaftMaschke2018,GernandtHallerReis2021,
MehrmannvanderSchaft2023,vanderSchaftMehrmann2023}; see also the recent results on nonlinear
systems in \cite{MehrmannSchaft2026}. In particular, \cite[Thm.~15]{GernandtPhilippPreusterSchaller2024}
established a geometric-to-descriptor conversion with internal and
external ports. More recently,
\cite{CaballeriaRamirezLefevreLeGorrec2026} developed equivalent
finite-dimensional representations of complete port-Hamiltonian
systems with Dirac and Lagrange constraints, including resistive
and external ports, and related them to constrained input--output
formulations; see also \cite{vanderSchaftMaschke2020}.

Algebraic constraints arising from Dirac structures are already
familiar in finite dimensions. The additional difficulty in infinite
dimensions is that the operator part of a Dirac relation may be
unbounded. Its domain can encode spatial regularity and boundary
conditions that cannot be recovered merely from the algebraic
constraints. A conversion to descriptor form must preserve these
domain restrictions as well. Moreover, the spaces associated with
algebraic constraints may themselves be infinite-dimensional, as in
models of incompressible continua.

Every closed linear relation admits bounded image and kernel
representations. The challenge here is therefore not simply to write
the geometric inclusions as operator equations. Rather, we seek a
structured descriptor representation in which the Dirac
interconnection is realized by a skew-adjoint operator on enlarged
Hilbert spaces. The transition from a structured descriptor equation
to geometric relations is comparatively direct: the inverse
graph of the interconnection operator gives a Dirac structure.
For the converse direction, we use the operator--multivalued
decomposition of the inverse Dirac relation to obtain a
skew-adjoint operator representation. An auxiliary variable
accounts for the multivalued part, while the original operator
domain is retained. The construction draws on the theory of linear relations
\cite{Arens1961,Cross1998,HassiSnooSzafraniec2009,
BehrndtHassiDeSnoo2020}; related infinite-dimensional Dirac
representations are studied in \cite{IftimeRomanSandovici2014}.
Combined with parametrizations of the storage and resistive relations,
this gives a closed descriptor operator with internal and
external ports. We prove an exact correspondence of classical trajectories,
including preservation of the power balance.

An incompressible Kelvin--Voigt solid illustrates how the
descriptor representation retains unbounded differential operators,
infinite-dimensional incompressibility constraints, and
operator-domain restrictions.

\section{Structure of skew-adjoint relations}
\label{sec:preliminaries}
We introduce notation for Hilbert spaces and linear relations,
recall the operator--multivalued decomposition of skew-adjoint
relations, and construct a skew-adjoint operator representation
on enlarged Hilbert spaces without additional domain regularity.

\subsection{Hilbert-space notation and duality}
All Hilbert spaces are complex, with inner products linear in the first
argument. The anti-dual $\U\adjun$ of a Hilbert space $\U$ consists of its
bounded antilinear functionals. We denote the canonical duality pairing
between $\U\adjun$ and $\U$ by
$\dualprod{u'}{u}_{\U\adjun,\U}$, and use
\[
 \dualprod{u}{u'}_{\U,\U\adjun}
 :=\overline{\dualprod{u'}{u}_{\U\adjun,\U}}
\]
for the reversed order. Like the inner products, both pairings are
linear in the first and antilinear in the second argument.
The Riesz isomorphism and the canonical biduality map are characterized by
\[
 \begin{aligned}
  \Rf_{\U}&:\U\to\U\adjun,
  &\dualprod{\Rf_{\U}u}{v}_{\U\adjun,\U}
  &=\scprod{u}{v}_{\U},\\
  \iota_{\U}&:\U\to\U\adjun[2],
  &\dualprod{\iota_{\U}u}{u'}_{\U\adjun[2],\U\adjun}
  &=\dualprod{u}{u'}_{\U,\U\adjun}.
 \end{aligned}
\]
Both maps are linear and unitary. For a bounded linear operator
$A:\U\to\W$, its anti-dual adjoint
$A\adjun:\W\adjun\to\U\adjun$ is characterized by
\[
 \dualprod{A\adjun w'}{u}_{\U\adjun,\U}
 =\dualprod{w'}{Au}_{\W\adjun,\W}.
\]
The same notation is used for densely defined operators, with
the anti-dual adjoint defined on its natural domain.

We use $\U\times\V$ for Cartesian products and
$\mathcal{S}_1\perp\mathcal{S}_2$ for orthogonal direct sums
of subspaces of Hilbert spaces. Hilbert-space products carry
the product inner product.

Let $\U,\V$ be Hilbert spaces and let $\Gamma:\V\to \U\adjun$ be a bounded
isomorphism. We use the duality pairing
$\dualprod{v}{u}_{\Gamma}
:=\dualprod{\Gamma v}{u}_{\U\adjun,\U}$ and denote the annihilator
of a subspace $\mathcal{S}\subseteq \U$ by
$\mathcal{S}^\circ:=\setdef*{v\in \V}{\dualprod{v}{s}_{\Gamma}=0\text{ for all }s\in \mathcal{S}}$.
By contrast, $\mathcal{S}^\perp$ denotes the Hilbert-space orthogonal complement.
Orthogonal complements and projections refer to the given inner products.

\subsection{Linear relations and their adjoints}
\label{sec:relations}
A linear relation from $\U$ to $\V$ is a subspace
$\mathcal T\subseteq \U\times \V$. We write $\dom\mathcal T$ and $\ran\mathcal T$ for its
domain and range. Its multivalued part is
$\setdef*{v\in\V}{(0,v)\in\mathcal T}$.
As usual, $\mathcal T^{-1}:=\setdef*{(v,u)}{(u,v)\in\mathcal T}$
and $-\mathcal T:=\setdef*{(u,-v)}{(u,v)\in\mathcal T}$.

For an operator $A:\dom A\subseteq \U\to \V$, we write
$\graph A:=\setdef*{(u,Au)}{u\in\dom A}$. A relation is an operator
graph exactly when its multivalued part vanishes.
For background on linear relations, see
\cite{Arens1961,Cross1998,HassiSebestyenDeSnooSzafraniec2007,
HassiSnooSzafraniec2009,BehrndtHassiDeSnoo2020}.

With respect to the duality map $\Gamma$, the adjoint relation of $\mathcal T$ is
\begin{equation}\label{eq:relation-adjoint}
 \mathcal T\adjunX{\Gamma}
 :=\setdef*{(a,c)\in \U\times \V}{
 \dualprod{c}{u}_{\Gamma}
 =\overline{\dualprod{v}{a}_{\Gamma}}
 \text{ for all }(u,v)\in\mathcal T}.
\end{equation}
We call $\mathcal T$ self-adjoint or skew-adjoint if
$\mathcal T\adjunX{\Gamma}=\mathcal T$ or
$\mathcal T\adjunX{\Gamma}=-\mathcal T$, respectively.
For a skew-adjoint relation, the multivalued part equals
$(\dom\mathcal T)^\circ=\mathcal{K}^\circ$, where
$\mathcal{K}:=\overline{\dom\mathcal T}$. Thus $\mathcal T$ is
single-valued if and only if its domain is dense in $\U$.

\subsection{Operator representations of skew-adjoint relations}
The possible multivaluedness of a skew-adjoint relation is thus
determined by the closure of its domain. In infinite dimensions,
$\dom\mathcal T$ may moreover be a proper dense subspace of its closure $\mathcal{K}$.
The following theorem treats these two features separately. Its first
part is the standard operator--multivalued decomposition of a
self-adjoint relation, transferred to the skew-adjoint setting
\cite{HassiSnooSzafraniec2009,BehrndtHassiDeSnoo2020}.
The second part realizes the entire relation by a skew-adjoint
operator without replacing its domain by its closure.

\begin{theorem}[Representation of skew-adjoint relations]\label{thm:canonical-lift}
Let $\mathcal T\subseteq \U\times \V$ be skew-adjoint with respect to
a unitary isomorphism $\Gamma:\V\to \U\adjun$. Set
$\mathcal{K}:=\overline{\dom\mathcal T}$ and $\mathcal{M}:=\mathcal{K}^\perp$.
Then the following statements hold.
\begin{enumerate}
\item\label{itm:skew-decomposition}
\emph{Operator--multivalued decomposition.}
There exists a unique closed operator
$J_{\mathcal{K}}:\dom\mathcal T\subseteq \mathcal{K}\to \mathcal{M}^\circ$ such that
\[
 \mathcal T=\graph J_{\mathcal{K}}\perp(\{0\}\times \mathcal{K}^\circ).
\]
Its zero extension $J_0$, with
$\dom J_0=\dom\mathcal T \perp\mathcal{M}$ and
$J_0(u_{\mathcal{K}}+u_\perp)=J_{\mathcal{K}}u_{\mathcal{K}}$ for $u_{\mathcal{K}}\in\dom\mathcal T$ and
$u_\perp\in \mathcal{M}$, is closed, densely defined, and skew-adjoint
with respect to $\Gamma$.

\item\label{itm:skew-lift}
\emph{Skew-adjoint operator representation.}
Equip $\Y:=\U/\mathcal{K}$ with its quotient Hilbert norm, let
$q:\U\to\Y$ be the quotient map, and set
$j:=\Gamma^{-1}q\adjun:\Y\adjun\to \V$. The block operator
\begin{equation}\label{eq:lift-block}
 \widehat J:=
 \begin{bmatrix}
  J_0&j\\
  -q&0
 \end{bmatrix},
 \qquad
 \dom\widehat J=\dom J_0\times\Y\adjun,
\end{equation}
from $\U\times\Y\adjun$ to $\V\times\Y$ is closed, densely
defined, and skew-adjoint with respect to the unitary duality map
$\diag(\Gamma,\iota_\Y)$ from $\V\times\Y$ to
$\U\adjun\times\Y\adjun[2]$. For every $(u,v)\in \U\times \V$,
\[
 (u,v)\in\mathcal T
 \quad\Longleftrightarrow\quad
 \exists!\,\lambda\in\Y\adjun:\quad
 \widehat J\begin{pmatrix}u\\\lambda\end{pmatrix}
 =\begin{pmatrix}v\\0\end{pmatrix}.
\]
Moreover, $j$ is an isometric isomorphism onto $\mathcal{K}^\circ$.
The unique auxiliary variable is $\lambda=j^{-1}P_{\mathcal{K}^\circ}v$
and satisfies $\norm{\lambda}_{\Y\adjun}\leq\norm{v}_{\V}$,
where $P_{\mathcal{K}^\circ}$ is the orthogonal projection in $\V$.
\end{enumerate}
\end{theorem}

\begin{proof}
\labelcref{itm:skew-decomposition} Set $W:=\Rf_{\U}^{-1}\Gamma:\V\to \U$. The relation
$\mathcal A:=\setdef*{(u,iWv)}{(u,v)\in\mathcal T}$ is self-adjoint
in $\U$ and has multivalued part $\mathcal{M}$. By
\cite[Thm.~1.5.1]{BehrndtHassiDeSnoo2020}, it decomposes as
$\mathcal A=\graph L_{\mathcal{K}}\perp(\{0\}\times \mathcal{M})$, where $L_{\mathcal{K}}$
is self-adjoint in $\mathcal{K}$ and $\dom L_{\mathcal{K}}=\dom\mathcal T$.
Since $W^{-1}\mathcal{K}=\mathcal{M}^\circ$ and $W^{-1}\mathcal{M}=\mathcal{K}^\circ$, the operator
$J_{\mathcal{K}}:=-iW^{-1}L_{\mathcal{K}}$ yields the asserted unique decomposition.
Moreover, $WJ_0$ acts as $-iL_{\mathcal{K}}$ on $\dom L_{\mathcal{K}}$ and vanishes on
$\mathcal{M}$. Since $\U=\mathcal{K}\perp \mathcal{M}$, this is a closed, densely defined,
skew-adjoint operator, and hence so is $J_0$ with respect to
$\Gamma$.

\labelcref{itm:skew-lift} The anti-dual adjoint $q\adjun$ maps
$\Y\adjun$ isometrically onto the annihilator of $\mathcal{K}$ in
$\U\adjun$. Thus $j=\Gamma^{-1}q\adjun$ maps
$\Y\adjun$ isometrically onto $\mathcal{K}^\circ$.
The off-diagonal operator
$N:=\sbmat{0}{j}{-q}{0}$, defined on
$\U\times\Y\adjun$, is bounded and skew-adjoint with respect to
$\diag(\Gamma,\iota_\Y)$, since $\Gamma j=q\adjun$.
The block-diagonal operator $\diag(J_0,0)$ is skew-adjoint as well.
By the bounded-perturbation adjoint formula,
$\widehat J=\diag(J_0,0)+N$ is therefore skew-adjoint on
$\dom J_0\times\Y\adjun$, hence closed and densely defined.

The second component of the lifted equation gives $qu=0$.
Together with $u\in\dom J_0$, this is equivalent to
$u\in\dom\mathcal T$; its first component is then
$v=J_{\mathcal{K}}u+j\lambda$. By \labelcref{itm:skew-decomposition}, this holds for a
unique $\lambda$ exactly when $(u,v)\in\mathcal T$.
Since $J_{\mathcal{K}}u\in \mathcal{M}^\circ$ and $j\lambda\in \mathcal{K}^\circ$,
orthogonal projection gives the multiplier formula and its
norm estimate.
\end{proof}

The equation $qu=0$ enforces the closed subspace
constraint $u\in\mathcal{K}$, whereas the domain of $J_0$
retains the possibly nonclosed domain of $\mathcal T$.
The construction depends on the chosen Hilbert-space
inner products through $\U=\mathcal{K}\perp\mathcal{K}^\perp$.

\begin{remark}[Converse of the decomposition]\label{rem:converse-decomposition}
Let $\mathcal{K}\subseteq \U$ be closed, and suppose that
$J_{\mathcal{K}}:\dom J_{\mathcal{K}}\subseteq \mathcal{K}\to(\mathcal{K}^\perp)^\circ$ has a zero
extension to $\dom J_{\mathcal{K}}\perp \mathcal{K}^\perp$ that is skew-adjoint
with respect to $\Gamma$. Then
$\graph J_{\mathcal{K}}\perp(\{0\}\times \mathcal{K}^\circ)$ is a skew-adjoint
relation. This follows by reversing the construction in
\labelcref{itm:skew-decomposition} of \Cref{thm:canonical-lift}.
\end{remark}

\begin{remark}[Finite-dimensional case]\label{rem:finite}
In finite dimensions, $\mathcal{K}=\dom\mathcal T$, and $J_0$ is everywhere
defined. Identify $\V$ with $\U$ via $W=\Rf_{\U}^{-1}\Gamma$ and
identify $\Y$ and $\Y\adjun$ with $\mathcal{M}=\mathcal{K}^\perp$ by the Hilbert
metrics. In these coordinates, $\widehat J$ becomes the
skew-Hermitian block matrix
\[
 \begin{bmatrix}
  WJ_0&\operatorname{incl}_{\mathcal{M}}\\
  -P_{\mathcal{M}}&0
 \end{bmatrix},
\]
where $\operatorname{incl}_{\mathcal{M}}:\mathcal{M}\to \U$ is the inclusion and
$P_{\mathcal{M}}:\U\to \mathcal{M}$ the orthogonal projection. In infinite dimensions,
$\dom\mathcal T$ can be strictly smaller than $\mathcal{K}$, so it cannot
in general be replaced by its closure.
\end{remark}

\section{Geometric port-Hamiltonian systems}
\label{sec:geometry}

\subsection{State and port spaces}
Let $\X$, $\Flow_{\mathrm R}$, and $\U$ be Hilbert spaces.
The storage flow space is $\Flow_{\mathrm S}:=\X$, with
$\Effort_{\mathrm S}:=\X\adjun$, so that the storage flow is
$\dot x$. We set $\Effort_{\mathrm R}:=\Flow_{\mathrm R}\adjun$
for the internal resistive port, while $u\in\U$ and
$y\in\U\adjun$ are the external input and output.
The storage and resistive efforts act by canonical duality;
the external port uses the reversed pairing
$\dualprod{u}{y}_{\U,\U\adjun}$ introduced above.
The total flow and effort spaces are the Cartesian products
\[
 \Flow_{\mathrm{tot}}
 :=\X\times\Flow_{\mathrm R}\times\U\adjun,\qquad
 \Effort_{\mathrm{tot}}
 :=\X\adjun\times\Effort_{\mathrm R}\times\U.
\]
With the canonical identification
$\Flow_{\mathrm{tot}}\adjun=
\X\adjun\times\Effort_{\mathrm R}\times\U\adjun[2]$,
the unitary duality map
$\Gamma_{\mathrm{tot}}:\Effort_{\mathrm{tot}}\to
\Flow_{\mathrm{tot}}\adjun$ is $\diag(I,I,\iota_{\U})$.
For total efforts $e$ and flows $f$, ordered by storage,
resistive, and external components, the corresponding pairing is
\[
 \dualprod{e}{f}_{\Gamma_{\mathrm{tot}}}
 =\dualprod{e_{\mathrm S}}{f_{\mathrm S}}_{\X\adjun,\X}
  +\dualprod{e_{\mathrm R}}{f_{\mathrm R}}_{\Effort_{\mathrm R},\Flow_{\mathrm R}}
  +\dualprod{u}{y}_{\U,\U\adjun}.
\]
For the reverse direction, we use the unitary map
$\Gamma_{\mathrm{tot}}^{\vee}:=
\Gamma_{\mathrm{tot}}\adjun\iota_{\Flow_{\mathrm{tot}}}$
from $\Flow_{\mathrm{tot}}$ to $\Effort_{\mathrm{tot}}\adjun$.
Its pairing satisfies
$\dualprod{f}{e}_{\Gamma_{\mathrm{tot}}^{\vee}}
=\overline{\dualprod{e}{f}_{\Gamma_{\mathrm{tot}}}}$.

\subsection{Dirac, Lagrangian, and resistive structures}
\label{subsec:geometric-phs}
The Dirac, Lagrangian, and resistive structures describe
power-conserving interconnection, state--effort relations,
and admissible internal port variables, respectively.

On $\Flow_{\mathrm{tot}}\times\Effort_{\mathrm{tot}}$ define the
Hermitian form
\[
 \sdprod{(f_1,e_1)}{(f_2,e_2)}_+
 :=\dualprod{e_1}{f_2}_{\Gamma_{\mathrm{tot}}}
   +\overline{\dualprod{e_2}{f_1}_{\Gamma_{\mathrm{tot}}}}.
\]
For a subspace $\mathcal S$ of this product space, let
$\mathcal S^{\perp_+}$ denote its orthogonal complement with respect
to this form. A \emph{Dirac structure} is a closed subspace
$\D\subseteq\Flow_{\mathrm{tot}}\times\Effort_{\mathrm{tot}}$
satisfying
\[
 \D=\D^{\perp_+}.
\]
By the definition of the relation adjoint in
\eqref{eq:relation-adjoint}, this is equivalent to
$\D\adjunX{\Gamma_{\mathrm{tot}}}=-\D$ when $\D$ is viewed
as a relation from flows to efforts. Thus the Dirac condition expresses maximality, beyond the
power identity $\Re\dualprod{e}{f}_{\Gamma_{\mathrm{tot}}}=0$
for $(f,e)\in\D$.

On $\X\times\X\adjun$ define the skew-Hermitian form
\[
 \sdprod{(x_1,e_1)}{(x_2,e_2)}_-
 :=\dualprod{x_1}{e_2}_{\X,\X\adjun}
   -\dualprod{e_1}{x_2}_{\X\adjun,\X}.
\]
For subspaces $\mathcal S\subseteq\X\times\X\adjun$,
write $\mathcal S^{\perp_-}$ for the corresponding
skew-orthogonal complement. A \emph{Lagrangian structure} is a
closed subspace $\Lagr\subseteq\X\times\X\adjun$
satisfying $\Lagr=\Lagr^{\perp_-}$, equivalently a
self-adjoint relation with respect to the canonical
state--effort pairing.

The Lagrangian structure determines the stored energy
$\mathcal H(x,e_{\mathrm S})
:=\tfrac12\dualprod{e_{\mathrm S}}{x}_{\X\adjun,\X}$ for
$(x,e_{\mathrm S})\in\Lagr$. By the Lagrangian property, this value is real and depends only
on $x$. Indeed, if
$(x,e_{\mathrm S}),(x,\widetilde e_{\mathrm S})\in\Lagr$,
then $(0,e_{\mathrm S}-\widetilde e_{\mathrm S})\in\Lagr$;
self-adjointness gives
$\dualprod{e_{\mathrm S}-\widetilde e_{\mathrm S}}{x}_{\X\adjun,\X}=0$.
We do not assume that this energy is positive.

The \emph{resistive relation} is a closed linear relation
$\Rc\subseteq\Flow_{\mathrm R}\times\Effort_{\mathrm R}$.
We call it \emph{dissipative} if
\begin{equation}\label{eq:R-resistive}
 \Re\dualprod{e_{\mathrm R}}{f_{\mathrm R}}_{
    \Effort_{\mathrm R},\Flow_{\mathrm R}}\leq 0
 \quad\text{for all }(f_{\mathrm R},e_{\mathrm R})\in\Rc.
\end{equation}
The conversion results below require only that $\Rc$ be closed;
dissipativity is relevant when interpreting the sign of its
contribution to the energy balance.

\subsection{Geometric systems and classical trajectories}
The Dirac structure couples the storage flow $\dot x$ to the
resistive port $(f_{\mathrm R},e_{\mathrm R})$ and the external
input $u$ and output $y$ without contributing net power.
Resistive and external flows are oriented opposite to the storage
flow, yielding the usual signs in the energy balance.

\begin{definition}[Geometric port-Hamiltonian system]
\label{def:geometric-system}
Let $\D\subseteq\Flow_{\mathrm{tot}}\times\Effort_{\mathrm{tot}}$
be a Dirac structure,
$\Lagr\subseteq\X\times\X\adjun$ a Lagrangian structure,
and $\Rc\subseteq\Flow_{\mathrm R}\times\Effort_{\mathrm R}$
a closed linear relation. The \emph{geometric port-Hamiltonian
system} associated with these relations is the differential-algebraic
inclusion \eqref{eq:geometric-system}, with variables in the spaces
specified above.
A \emph{classical trajectory} on an interval
$I\subseteq\R$ consists of functions
\[
\begin{aligned}
 x&\in\conC^1(I;\X),&
 f_{\mathrm R}&\in\conC(I;\Flow_{\mathrm R}),&
 u&\in\conC(I;\U),\\
 e_{\mathrm S}&\in\conC(I;\X\adjun),&
 e_{\mathrm R}&\in\conC(I;\Effort_{\mathrm R}),&
 y&\in\conC(I;\U\adjun).
\end{aligned}
\]
satisfying \eqref{eq:geometric-system} for every $t\in I$.
\end{definition}

The formulation is behavioral: the inclusion specifies the
admissible trajectories without requiring the state derivative or
the output to be uniquely determined by the other variables.
To explain the signs, Dirac isotropy implies that the storage power
equals the sum of resistive and external port powers. The Lagrangian
property gives the chain rule
$\ddts\mathcal H(x(t),e_{\mathrm S}(t))
=\Re\dualprod{e_{\mathrm S}(t)}{\dot x(t)}_{\X\adjun,\X}$.
Hence every classical trajectory
satisfies
\[
 \ddts\mathcal H(x(t),e_{\mathrm S}(t))
 =\Re\dualprod{e_{\mathrm R}(t)}{f_{\mathrm R}(t)}_{\Effort_{\mathrm R},\Flow_{\mathrm R}}
  +\Re\dualprod{u(t)}{y(t)}_{\U,\U\adjun}.
\]
The resistive contribution is nonpositive when $\Rc$ is dissipative.
The storage chain rule is justified under the stated regularity
assumptions in \Cref{sec:ph-structure}.

\section{Between geometric systems and descriptor equations}
\label{sec:parametrization}
To convert between \eqref{eq:geometric-system} and
\eqref{eq:compact-descriptor}, we first parametrize the storage
and resistive relations, then treat the Dirac interconnection.

\subsection{Parametrizing storage and resistive relations}
\label{sec:params}

Since $\Lagr$ and $\Rc$ are closed linear relations,
there exist Hilbert spaces $\Z_{\Lagr}$ and $\Z_{\mathcal R}$
and bounded isomorphisms
$(E_{\Lagr},Q_{\Lagr}):\Z_{\Lagr}\to\Lagr$ and
$(F_{\mathcal R},E_{\mathcal R}):\Z_{\mathcal R}\to\Rc$ \cite[Thm.~1.10.1 and Prop.~1.10.3]{BehrndtHassiDeSnoo2020}.
Thus $(x,e_{\mathrm S})=(E_{\Lagr}z,Q_{\Lagr}z)$ and
$(f_{\mathrm R},e_{\mathrm R})=(F_{\mathcal R}r,E_{\mathcal R}r)$.
For $\Lagr=\graph A$ with $A$ closed and possibly unbounded,
one may take $\Z_{\Lagr}=\dom A$ with its graph norm and
$(E_{\Lagr},Q_{\Lagr})=(I,A)$.
Only $E_{\Lagr}z$, not necessarily $z$ in the graph norm,
is differentiated in the descriptor equation.

The Dirac structure $\D$ also admits bounded image
parametrizations. Here we seek a \emph{skew-adjoint operator}
representation of the interconnection on enlarged effort and
flow spaces, using \Cref{thm:canonical-lift}.

\subsection{From operator-based descriptors to geometric systems}
\label{subsec:descriptor-to-geometric}
\begin{proposition}[Geometric interpretation of a structured descriptor equation]
\label{prop:descriptor-to-geometric}
Let $J:\dom J\subseteq\Effort_{\mathrm{tot}}\to\Flow_{\mathrm{tot}}$
be densely defined and skew-adjoint with respect to
$\Gamma_{\mathrm{tot}}^{\vee}$. Let $(E_{\Lagr},Q_{\Lagr})$ and
$(F_{\mathcal R},E_{\mathcal R})$ be bounded isomorphisms
onto a Lagrangian structure $\Lagr$ and a closed resistive
relation $\Rc$, respectively. Consider
\[
  \begin{pmatrix}
  \ddts(E_{\Lagr}z(t))\\-F_{\mathcal R}r(t)\\-y(t)
 \end{pmatrix}
 =Je(t),\quad
 e(t)=\begin{pmatrix}
  Q_{\Lagr}z(t)\\E_{\mathcal R}r(t)\\u(t)
 \end{pmatrix}\in\dom J.
\]
Then $\D:=\setdef*{(Je,e)}{e\in\dom J}$ is a Dirac structure.
Under the regularity assumptions
$E_{\Lagr}z\in\conC^1(I;\X)$ and continuity of
$z,r,u,y$ in their respective spaces, solutions of the
displayed equation correspond bijectively to classical
trajectories of \eqref{eq:geometric-system}, with
$x=E_{\Lagr}z$, $e_{\mathrm S}=Q_{\Lagr}z$, $f_{\mathrm R}=F_{\mathcal R}r$,
and $e_{\mathrm R}=E_{\mathcal R}r$.
\end{proposition}
\begin{proof}
The inverse of the skew-adjoint graph of $J$ is a Dirac
structure. Under the two bijective parametrizations, the displayed
equation is precisely the Dirac inclusion in
\eqref{eq:geometric-system}; the other inclusions hold by
construction.
\end{proof}

\subsection{From geometric Dirac structures to operators}
\label{subsec:geometric-to-descriptor}
For a general Dirac structure $\D$, the inverse relation
$\mathcal T_{\D}:=\D^{-1}\subseteq
\Effort_{\mathrm{tot}}\times\Flow_{\mathrm{tot}}$
need not be single-valued. With the reversed duality
$\Gamma_{\mathrm{tot}}^{\vee}$ from \Cref{sec:geometry}, we have
$(\D^{-1})\adjunX{\Gamma_{\mathrm{tot}}^{\vee}}
=(\D\adjunX{\Gamma_{\mathrm{tot}}})^{-1}
=-\D^{-1}$; see also
\cite[Prop.~2.8(i)]{IftimeRomanSandovici2014}.
Thus $\mathcal T_{\D}$ is skew-adjoint and
\Cref{thm:canonical-lift} applies. Put
$\mathcal{K}_{\D}:=\overline{\dom\mathcal T_{\D}}
=\overline{\ran\D}$. Its multivalued part is
$\mathcal{K}_{\D}^{\circ}$ by \Cref{thm:canonical-lift}.

\begin{corollary}[Operator representation of the Dirac structure]
\label{thm:dirac-operator-representation}
Let $\Y_{\D}:=\Effort_{\mathrm{tot}}/\mathcal{K}_{\D}$
with quotient map $q_{\D}:\Effort_{\mathrm{tot}}\to\Y_{\D}$,
and set $j_{\D}:=(\Gamma_{\mathrm{tot}}^{\vee})^{-1}q_{\D}\adjun$.
Applying \eqref{eq:lift-block} to $\mathcal T_{\D}$ with
$\Gamma_{\mathrm{tot}}^{\vee}$ gives a closed, densely defined
operator $\widehat J_{\D}$ with domain
$\dom\widehat J_{\D}
=(\dom\mathcal T_{\D}\perp\mathcal{K}_{\D}^{\perp})
\times\Y_{\D}\adjun$,
skew-adjoint with respect to
$\diag(\Gamma_{\mathrm{tot}}^{\vee},\iota_{\Y_{\D}})$.
Moreover, $(f,e)\in\D$ if and only if there is a unique
$\lambda\in\Y_{\D}\adjun$ satisfying the first equation below;
the second gives its recovery:
\begin{gather}
 \widehat J_{\D}\pvek{e}{\lambda}=\pvek{f}{0},
 \label{eq:D-recovery}\\
 \lambda=\Lambda_{\D}(f,e)
 :=j_{\D}^{-1}P_{\mathcal{K}_{\D}^{\circ}}f,
 \label{eq:Lambda-D}
\end{gather}
where $P_{\mathcal{K}_{\D}^{\circ}}$ is the orthogonal projection in
$\Flow_{\mathrm{tot}}$. The recovery map is bounded and satisfies
$\norm{\Lambda_{\D}(f,e)}_{\Y_{\D}\adjun}
\le\norm{f}_{\Flow_{\mathrm{tot}}}$.
\end{corollary}
\begin{proof}
Apply \Cref{thm:canonical-lift} to
$\mathcal T=\D^{-1}\subseteq\Effort_{\mathrm{tot}}\times\Flow_{\mathrm{tot}}$
with $\Gamma=\Gamma_{\mathrm{tot}}^{\vee}$.
\end{proof}

\begin{remark}[Effort constraints and operator domains]
\label{rem:dense-vs-algebraic}
The quotient equation $q_{\D}e=0$ in \eqref{eq:D-recovery}
enforces $e\in\mathcal{K}_{\D}$. Together with the domain of
$\widehat J_{\D}$, it also enforces $e\in\dom\mathcal T_{\D}$, since
$(\dom\mathcal T_{\D}\perp\mathcal{K}_{\D}^{\perp})
\cap\mathcal{K}_{\D}=\dom\mathcal T_{\D}$.
Thus the operator domain retains spatial regularity and boundary
conditions independently of algebraic constraints. If $\mathcal{K}_{\D}=\Effort_{\mathrm{tot}}$, no multiplier
is needed ($\Y_{\D}=\{0\}$), even when
$\dom\mathcal T_{\D}$ is a proper dense subspace.
Moreover, $\widehat J_{\D}$ is bounded precisely when
$\dom\mathcal T_{\D}=\mathcal{K}_{\D}$.
\end{remark}

\subsection{Construction and exact correspondence}
\label{sec:general-descriptor}
Combining the parametrizations from \Cref{sec:params} with
\eqref{eq:D-recovery}, define the bounded map
$S:\Effort_{\mathrm S}\times\Z_{\mathcal R}
\times\Y_{\D}\adjun\times\U\to
\Effort_{\mathrm{tot}}\times\Y_{\D}\adjun$ by
\[
 S\begin{pmatrix}e\\r\\\lambda\\u\end{pmatrix}
 =\begin{pmatrix}e\\E_{\mathcal R}r\\u\\\lambda\end{pmatrix}.
\]
Let $\Pi$ exchange the last two components of
$\Flow_{\mathrm{tot}}\times\Y_{\D}$.
For $w$ with components $e,r,\lambda,u$, define
\begin{align}
 \dom M&:=S^{-1}(\dom\widehat J_{\D}),
   \label{eq:M-domain}\\
 Mw&:=\Pi\widehat J_{\D}Sw+
 \begin{pmatrix}0\\F_{\mathcal R}r\\0\\0\end{pmatrix}
 \quad(w\in\dom M).
   \label{eq:M-def}
\end{align}
The zero rows impose the resistive and Dirac constraints.
The finite-dimensional construction with internal and external
ports is given in \cite[Thm.~15]{GernandtPhilippPreusterSchaller2024}.

\begin{theorem}[Equivalence of the geometric and descriptor formulations]
\label{thm:geometric-compact-equivalence}
Let $\D$ be a Dirac structure, $\Lagr$ a Lagrangian
structure, and $\Rc$ a closed resistive relation, not necessarily
dissipative. Choose bounded isomorphisms $(E_{\Lagr},Q_{\Lagr})$ and
$(F_{\mathcal R},E_{\mathcal R})$ as in \Cref{sec:params}
and construct $\widehat J_{\D}$ and $M$ as above.
Then $M$ is closed.

Set $\Z:=\Z_{\Lagr}\times\Z_{\mathcal R}\times\Y_{\D}\adjun$.
Define bounded operators $E:\Z\to\Flow_{\mathrm S}\times\Flow_{\mathrm R}\times\Y_{\D}$
and $Q:\Z\to\Effort_{\mathrm S}\times\Z_{\mathcal R}\times\Y_{\D}\adjun$ by
\[
 E\zeta=\begin{pmatrix}E_{\Lagr}z\\0\\0\end{pmatrix},\quad
 Q\zeta=\begin{pmatrix}Q_{\Lagr}z\\r\\\lambda\end{pmatrix},\quad
 \text{where }\,\zeta=\begin{pmatrix}z\\r\\\lambda\end{pmatrix}.
\]
Then $x=E_{\Lagr}z$, $e_{\mathrm S}=Q_{\Lagr}z$,
$f_{\mathrm R}=F_{\mathcal R}r$, and
$e_{\mathrm R}=E_{\mathcal R}r$ define a bijection
between classical geometric trajectories on $I$ and triples
$\zeta\in\conC(I;\Z)$, $u\in\conC(I;\U)$,
$y\in\conC(I;\U\adjun)$ such that
$E_{\Lagr}z\in\conC^1(I;\X)$,
$\spvek{Q\zeta(t)}{u(t)}\in\dom M$,
and \eqref{eq:compact-descriptor} holds at every $t\in I$.
We call such triples \emph{classical descriptor trajectories}.
The external input and output are preserved.
\end{theorem}
\begin{proof}
\emph{Closedness.} Suppose $w_k\in\dom M$,
$w_k\to w$, and $Mw_k\to g$. Since $S$ and
$F_{\mathcal R}$ are bounded, $Sw_k\to Sw$,
and \eqref{eq:M-def} implies that
$\widehat J_{\D}Sw_k$ converges.
Closedness of $\widehat J_{\D}$ gives
$w\in\dom M$, and \eqref{eq:M-def} yields $Mw=g$.

\emph{Correspondence.} For a classical descriptor trajectory,
the parametrizations yield the regularity and membership
in $\Lagr$ and $\Rc$. By \eqref{eq:M-domain}--\eqref{eq:M-def} and
\eqref{eq:D-recovery}, the descriptor equation reproduces the
Dirac inclusion in \eqref{eq:geometric-system}.
Conversely, a classical
geometric trajectory determines unique continuous $z$ and $r$
by the bounded inverse parametrizations. The bounded map
\eqref{eq:Lambda-D} supplies a unique continuous
$\lambda=\Lambda_{\D}(f,e)$. These variables satisfy the
descriptor equation and its domain condition. The uniqueness
of $z,r,\lambda$ proves injectivity.
\end{proof}

\begin{remark}[Non-dense descriptor domains]
Although $\widehat J_{\D}$ is densely defined,
$\dom M=S^{-1}(\dom\widehat J_{\D})$ need not be dense.
\end{remark}

\begin{remark}[Port-free form]\label{rem:portfree}
Without internal and external ports, the construction reduces
to $\ddts(E\zeta(t))=\widehat J_{\D}Q\zeta(t)$, where
$E\spvek{z}{\lambda}=\spvek{E_{\Lagr}z}{0}$ and
$Q\spvek{z}{\lambda}=\spvek{Q_{\Lagr}z}{\lambda}$.
The interconnection operator remains closed and skew-adjoint.
\end{remark}

\subsection{Power structure and balance}
\label{sec:ph-structure}
The descriptor representation retains the geometric power
structure. Let $w$ have components $e,r,\lambda,u$, and let
$g_1,\ldots,g_4$ denote the storage, resistive, constraint,
and external components of a descriptor flow $g$.
The bounded sesquilinear effort--flow pairing is defined by
\begin{equation}\label{eq:power-pairing}
  \dualprod{w}{g}_{\Theta}
 :=\dualprod{e}{g_1}_{\X\adjun,\X}
   +\dualprod{E_{\mathcal R}r}{g_2}_{\Effort_{\mathrm R},\Flow_{\mathrm R}}
   +\dualprod{\lambda}{g_3}_{\Y_{\D}\adjun,\Y_{\D}}
   +\dualprod{u}{g_4}_{\U,\U\adjun}.
\end{equation}
This pairing may be degenerate.

\begin{proposition}[Port-Hamiltonian structure]
\label{prop:ph-structure}
Under the assumptions of
\Cref{thm:geometric-compact-equivalence}:
\begin{enumerate}
\item\label{itm:ph-dissipation}
\emph{Dissipation identity.}
For every $w\in\dom M$,
\[
 \Re\dualprod{w}{Mw}_{\Theta}
 =\Re\dualprod{E_{\mathcal R}r}{F_{\mathcal R}r}
 _{\Effort_{\mathrm R},\Flow_{\mathrm R}}.
\]
The resistive contribution is nonpositive for all $r$
if and only if $\Rc$ is dissipative.
\item\label{itm:ph-balance}
\emph{Power balance.}
Every classical descriptor trajectory satisfies
\begin{equation}\label{eq:classical-power}
 \ddts\mathcal H(E_{\Lagr}z(t),Q_{\Lagr}z(t))
 =\Re\dualprod{E_{\mathcal R}r(t)}{F_{\mathcal R}r(t)}
 _{\Effort_{\mathrm R},\Flow_{\mathrm R}}
 +\Re\dualprod{u(t)}{y(t)}_{\U,\U\adjun},
\end{equation}
where $t\mapsto\mathcal H(E_{\Lagr}z(t),Q_{\Lagr}z(t))$ is continuously differentiable.
\end{enumerate}
\end{proposition}
\begin{proof}
For \labelcref{itm:ph-dissipation}, the pairing of $w$
with $\Pi\widehat J_{\D}Sw$ is the conjugate of
the pairing of $Sw$ with $\widehat J_{\D}Sw$,
whose real part vanishes by skew-adjointness.
The second term in \eqref{eq:M-def} contributes
$\dualprod{E_{\mathcal R}r}{F_{\mathcal R}r}
_{\Effort_{\mathrm R},\Flow_{\mathrm R}}$
by \eqref{eq:power-pairing}.
The identity follows, and the sign characterization
comes from the bijective parametrization of $\Rc$
and \eqref{eq:R-resistive}.

For \labelcref{itm:ph-balance}, write $x=E_{\Lagr}z$ and $e_{\mathrm S}=Q_{\Lagr}z$.
The self-adjointness of $\Lagr$ gives
$\ddts\mathcal H(x,e_{\mathrm S})
=\Re\dualprod{e_{\mathrm S}}{\dot x}_{\X\adjun,\X}$ without differentiating
the effort. Indeed, for $s,t\in I$,
\[
 2\left(\mathcal H(x(t),e_{\mathrm S}(t))
  -\mathcal H(x(s),e_{\mathrm S}(s))\right)
    =\Re\dualprod{e_{\mathrm S}(t)+e_{\mathrm S}(s)}
     {x(t)-x(s)}_{\X\adjun,\X}.
\]
Difference quotients justify the chain rule under the stated
regularity. Finally, the descriptor equation and
\labelcref{itm:ph-dissipation} yield \eqref{eq:classical-power}.
\end{proof}

\section{Example: An incompressible Kelvin--Voigt solid with a boundary port}
\label{sec:kelvin}

We consider an incompressible Kelvin--Voigt solid with a
traction boundary port, illustrating unbounded Dirac operators,
infinite-dimensional constraints, and dissipation; see also
\cite[Sec.~3.2]{BrugnoliHaineMatignon2023}.

\paragraph{The mechanical model}
Let $\Omega\subset\R^d$, $d\in\{2,3\}$, be a bounded,
connected Lipschitz domain. Its boundary is decomposed,
up to sets of surface measure zero, into disjoint relatively
open parts $\Gamma_D$ and $\Gamma_P$, each of positive
surface measure.

Let $\rho,\mu,\eta\in\Lp{\infty}(\Omega)$ denote the mass
density, shear modulus, and shear viscosity, respectively.
These coefficients are real-valued and satisfy
$\rho\ge\rho_->0$, $\mu\ge\mu_->0$, and $\eta\ge0$
almost everywhere.
We denote the displacement by $q$, the velocity by
$v=\partial_tq$, and the momentum density by $p=\rho v$.
The linearized strain tensor is defined by
$\varepsilon(w):=\frac12(\nabla w+\nabla w^\top)$.
The linearized incompressible Kelvin--Voigt equations are
\begin{equation}\label{kv:eq:pde}
 \begin{aligned}
  \partial_tq(t)&=v(t), &
  \partial_tp(t)&=\operatorname{div}T(t),\\
  p(t)&=\rho v(t),&
  \operatorname{div}v(t)&=0,\\
  T(t)&=2\mu\varepsilon(q(t))
       +2\eta\varepsilon(v(t))-\pi(t)I.
 \end{aligned}
\end{equation}
Here $\pi$ is the pressure multiplier enforcing
incompressibility, and $T$ is the total stress tensor,
comprising elastic, viscous, and pressure contributions.

The boundary conditions and external output are
\begin{equation}\label{kv:eq:boundary}
 q(t)|_{\Gamma_D}=0,\qquad
 T(t)n|_{\Gamma_P}=u(t),\qquad
 y(t)=\gamma_Pv(t).
\end{equation}
Here $n$ is the outward unit normal, $\gamma_P$ the trace
on $\Gamma_P$, $u$ the prescribed traction, and $y$ the
boundary velocity. The traction condition is understood weakly;
without body forces, all external power enters through $\Gamma_P$.

\paragraph{Velocity, trace, and constraint spaces}
Let $\gamma_D$ denote the trace on $\Gamma_D$. Put
\[
 \begin{aligned}
  H&:=\Lp{2}(\Omega;\C^d),&
  \Hsym&:=\Lp{2}(\Omega;\C_{\rm sym}^{d\times d}),\\
  W&:=\setdef*{w\in\Hk{1}(\Omega;\C^d)}
                  {\gamma_Dw=0},&
  V&:=\setdef*{w\in W}{\operatorname{div}w=0}.
 \end{aligned}
\]
Let $G:=\gamma_PW$, equipped with the quotient norm
\[
 \norm{g}_G:=
 \inf_{\substack{w\in W\\\gamma_Pw=g}}
       \norm{w}_{\Hk{1}(\Omega;\C^d)}.
\]
Equivalently, $G$ consists of traces on $\Gamma_P$
whose extensions by zero to $\Gamma_D$ belong to
$\Hk{1/2}(\partial\Omega;\C^d)$.
Define the subspace of traces with zero total flux by
\[
 G_\sigma:=\setdef*{g\in G}
 {\int_{\Gamma_P}g\cdot n\,\mathrm dS=0}.
\]
The divergence theorem gives
$\gamma_PV\subseteq G_\sigma$.
Conversely, any $g\in G_\sigma$ can be extended
to $w\in W$, whose divergence has zero mean.
Using Bogovskii's operator
\cite[Thm.~III.3.1]{Galdi2011},
this divergence can be corrected without changing
the boundary trace, yielding a field in $V$ with
trace $g$. Hence $\gamma_PV=G_\sigma$.
Since the total-flux functional
$g\mapsto\int_{\Gamma_P}g\cdot n\,\mathrm dS$
is nonzero on $G$, its kernel $G_\sigma$
has codimension one in $G$.

We deliberately choose the full space $G$, rather than
$G_\sigma$, for the external effort and flow pairing
$\U:=G\adjun$, $\U\adjun\cong G$.
Under the canonical identification $G\cong\U\adjun$,
the two pairings agree via
$\dualprod{u}{g}_{\U,G}
 =\dualprod{u}{\iota_Gg}_{\U,\U\adjun}$ for all
$u\in\U,\ g\in G$.
Here $u\in\U$ is the boundary traction, whereas
$y=\gamma_Pv\in G_\sigma$ is identified with an element of
$\U\adjun$. Using $G$ rather than $G_\sigma$ retains the
constant-normal traction component and thus fixes the pressure level.

Define the solenoidal subspace and its orthogonal complement by
\[
 H_\sigma:=\overline V^{\,H},\qquad
 H_{\nabla,P}:=\setdef*{\nabla\varphi}
 {\varphi\in\Hk{1}(\Omega),\ \gamma_P\varphi=0}.
\]
The mixed-boundary Helmholtz decomposition reads
\begin{equation}\label{kv:eq:helmholtz}
 H=H_\sigma\perp H_{\nabla,P}.
\end{equation}
The divergence map
$\operatorname{div}:W\to\Lp{2}(\Omega)$ is
surjective: Bogovskii's operator provides all mean-zero
divergences, while a field in $W$ with nonzero flux
through $\Gamma_P$ supplies the constant component.
Since $V=\ker(\operatorname{div}|_W)$,
every $h\in V^\perp$ induces a well-defined
functional on $\Lp{2}(\Omega)$ through
$\operatorname{div}w\mapsto\scprod{h}{w}_H$.
By the open mapping and Riesz representation
theorems, this functional is represented by
some $\varphi\in\Lp{2}(\Omega)$.
Integration by parts yields
$h=-\nabla\varphi$, with
$\varphi\in\Hk{1}(\Omega)$ and
$\gamma_P\varphi=0$.
The converse follows from the Green identity,
proving \eqref{kv:eq:helmholtz}.
The space $H_{\nabla,P}$ is infinite-dimensional.

By Korn's first inequality for fields with vanishing trace
on a nonempty Dirichlet boundary portion
\cite[Cor.~8.8]{NeffPaulyWitsch2014}, the symmetric-gradient operator
\[
 A:V\subset H_\sigma\longrightarrow\Hsym,\qquad
 Aw:=\varepsilon(w),
\]
is closed and densely defined. The trace $\gamma_P:V\to G_\sigma$ is unbounded
for the $H$-norm: otherwise the compact embedding
$V\hookrightarrow H_\sigma$ would make this surjective
trace map compact, contradicting $\dim G_\sigma=\infty$.

\paragraph{A boundary-coupled Dirac operator}
To include the traction and both elastic and resistive stresses
in a single adjoint, introduce
\[
 C:V\subset H_\sigma\longrightarrow
 \Hsym\times\Hsym\times G,\qquad
  Cw:=\begin{pmatrix}Aw\\Aw\\-\gamma_Pw\end{pmatrix}.
\]
This operator is densely defined and closed: convergence in its
graph norm implies convergence in $W$ by Korn's inequality.
We identify the $L^2$-based spaces $H$ and $\Hsym$
with their anti-duals through the Riesz isomorphisms.
We also identify $G$ with $\U\adjun$ by canonical
biduality. With these identifications, the anti-dual adjoint $C\adjun$
is characterized by
\begin{equation}\label{kv:eq:adjoint}
 \scprod{C\adjun\begin{pmatrix}\sigma_s\\\sigma_R\\u\end{pmatrix}}{w}_H
 =\scprod{\sigma_s+\sigma_R}{Aw}_{\Hsym}
 -\dualprod{u}{\gamma_Pw}_{\U,G}
 \qquad(w\in V).
\end{equation}
The displayed column belongs to $\dom C\adjun$ precisely when
the right-hand side is bounded in the $H$-norm on $V$;
\eqref{kv:eq:adjoint} then determines its image uniquely.
In particular, the domain condition couples the \emph{sum} of
the elastic and viscous stresses with the boundary traction.
Neither stress need have a normal trace by itself.

For the geometric system choose the storage and resistive spaces
and relations
\[
 \begin{aligned}
  \X&:=\Hsym\times H,
    &\Lagr&:=\graph\diag(2\mu,\rho^{-1}),\\
  \Flow_{\mathrm R}&=\Effort_{\mathrm R}:=\Hsym,
    &\Rc&:=\setdef*{(-r,2\eta r)}{r\in\Hsym}.
 \end{aligned}
\]
The storage state is $x=\spvek{s}{p}$, with effort
$e_{\mathrm S}=\spvek{2\mu s}{\rho^{-1}p}$, where $s$ denotes
the strain variable and $p$ the momentum density. The energy is
\[
 \mathcal H(s,p)=\tfrac12\scprod{2\mu s}{s}_{\Hsym}
                 +\tfrac12\scprod{\rho^{-1}p}{p}_H.
\]
Both multipliers are bounded, and $\Rc$ is a closed dissipative
relation because $\eta\ge0$.

The interconnection will restrict the velocity effort
$e_p=\rho^{-1}p$ to $V$. Since $\overline V^{\,H}=H_\sigma$
and $\dom C\adjun$ is dense in
$\Hsym\times\Hsym\times\U$, the closed effort subspace
and its annihilator are
\[
 \mathcal K:=\Hsym\times H_\sigma\times\Hsym\times\U,
 \qquad
 \mathcal K^\circ
 :=\{0\}\times H_{\nabla,P}\times\{0\}\times\{0\}.
\]
The expression for $\mathcal K^\circ$ follows from
\eqref{kv:eq:helmholtz}. On $\mathcal K$ define
\[
 \dom J_{\mathcal K}
 :=\setdef*{e\in\mathcal K}{
    e_p\in V,\;
    \begin{pmatrix}e_s\\e_{\mathrm R}\\u\end{pmatrix}\in\dom C\adjun
   },\quad
 J_{\mathcal K}\begin{pmatrix}e_s\\e_p\\e_{\mathrm R}\\u\end{pmatrix}
 :=\begin{pmatrix}
    Ae_p\\
    -C\adjun\left(\begin{smallmatrix}e_s\\e_{\mathrm R}\\u\end{smallmatrix}\right)\\
    Ae_p\\-\gamma_Pe_p
   \end{pmatrix}.
\]
The last component is understood in $\U\adjun$.
After placing the velocity effort first and grouping the two
stress efforts with the boundary traction, and reordering the
flow coordinates accordingly, $J_{\mathcal K}$ takes the block form
$\left[\begin{smallmatrix}
  0&-C\adjun\\
  C&0
 \end{smallmatrix}
\right]$
on $\dom C\times\dom C\adjun$.
Since $C$ is closed and densely defined, its bidual
adjoint coincides with $C$ under the canonical
identifications. Thus the block operator is skew-adjoint
with respect to the stated duality. Its zero extension to
$\dom J_{\mathcal K}\perp\mathcal K^\perp$ is skew-adjoint as well.
Consequently,
\[
 \D:=\setdef*{(J_{\mathcal K}e+g,e)}
 {e\in\dom J_{\mathcal K},\ g\in\mathcal K^\circ}
\]
is a Dirac structure by \Cref{rem:converse-decomposition}.
The two restrictions
$\mathcal K\subsetneq\Effort_{\mathrm{tot}}$ and
$\dom J_{\mathcal K}\subsetneq\mathcal K$
exhibit the algebraic incompressibility constraint
and the unbounded spatial and boundary-port
operator domain, respectively.

\paragraph{Descriptor representation}
By the Helmholtz decomposition \eqref{kv:eq:helmholtz}, the quotient space
$\Y_{\D}=\Effort_{\mathrm{tot}}/\mathcal K$
can be identified with $H_{\nabla,P}$. We identify
$\Y_{\D}\adjun$ with the same space by Riesz and write
$\lambda\in H_{\nabla,P}$. Let $P_\sigma\in\Lb(H)$
denote the orthogonal projection onto $H_\sigma$ from
\eqref{kv:eq:helmholtz}. The maps $q_{\D}$ and $j_{\D}$ are
\[
 q_{\D}e=(I-P_\sigma)e_p,\qquad
 j_{\D}\lambda=\begin{pmatrix}0\\\lambda\\0\\0\end{pmatrix}.
\]
Using the storage parametrization
$E_{\Lagr}=I$, $Q_{\Lagr}=\diag(2\mu,\rho^{-1})$
and the resistive parametrization
$F_{\mathcal R}r=-r$, $E_{\mathcal R}r=2\eta r$,
the general construction \eqref{eq:M-domain}--\eqref{eq:M-def}
gives the descriptor system
\begin{equation}\label{kv:eq:descriptor}
  \begin{pmatrix}
   \partial_ts(t)\\\partial_tp(t)\\0\\0\\-y(t)
  \end{pmatrix}
  =M\begin{pmatrix}
    2\mu s(t)\\\rho^{-1}p(t)\\r(t)\\\lambda(t)\\u(t)
   \end{pmatrix},\qquad
 M\begin{pmatrix}e_s\\e_p\\r\\\lambda\\u\end{pmatrix}
 =\begin{pmatrix}
  AP_\sigma e_p\\
  -C\adjun\left(\begin{smallmatrix}e_s\\2\eta r\\u\end{smallmatrix}\right)+\lambda\\
  AP_\sigma e_p-r\\
  -(I-P_\sigma)e_p\\
  -\gamma_PP_\sigma e_p
 \end{pmatrix}.
\end{equation}
With $w$ denoting the effort column in \eqref{kv:eq:descriptor},
the domain of $M$ on
$\mathcal{W}:=\Hsym\times H\times\Hsym\times H_{\nabla,P}\times\U$ is
\[
  \dom M:=\setdef*{w\in\mathcal{W}}{
      P_\sigma e_p\in V,\;
      \begin{pmatrix}e_s\\2\eta r\\u\end{pmatrix}\in\dom C\adjun
   }.
\]
The third and fourth rows of \eqref{kv:eq:descriptor} enforce
$r=Av$ and $v=P_\sigma v\in V$, respectively, with
$v=\rho^{-1}p$. The fifth row yields $y=\gamma_Pv$.
Thus the input is a boundary traction and the output a
velocity trace, not a distributed body-force/velocity pair.

\paragraph{Equivalence and power balance}
We first show that classical descriptor trajectories of
\eqref{kv:eq:descriptor} with
compatible initial strain recover \eqref{kv:eq:pde} and its
boundary conditions.
Set $v=\rho^{-1}p$ and
$\sigma:=2\mu s+2\eta r$.
Since $\lambda\in H_{\nabla,P}=V^\perp$,
the second row of \eqref{kv:eq:descriptor}
and \eqref{kv:eq:adjoint} yield
\[
  \scprod{\partial_tp(t)}{w}_H
 +\scprod{\sigma(t)}{\varepsilon(w)}_{\Hsym}
 -\dualprod{u(t)}{\gamma_Pw}_{\U,G}=0
 \qquad(w\in V).
\]

To recover the pressure, observe that the
left-hand side defines a bounded functional
on $W$ that vanishes on
$V=\ker(\operatorname{div}|_W)$.
Since $\operatorname{div}:W\to\Lp{2}(\Omega)$
is bounded and surjective, the open mapping
and Riesz representation theorems yield, for each
$t\in I$, a unique $\pi(t)\in\Lp{2}(\Omega)$ such that
\begin{equation}\label{kv:eq:weakpressure}
  \scprod{\partial_tp(t)}{w}_H
 +\scprod{\sigma(t)}{\varepsilon(w)}_{\Hsym}
 -\dualprod{u(t)}{\gamma_Pw}_{\U,G}
 =\scprod{\pi(t)}{\operatorname{div}w}_{\Lp{2}(\Omega)}
 \qquad(w\in W).
\end{equation}
Hence $T(t):=\sigma(t)-\pi(t)I$ satisfies
$\operatorname{div}T(t)=\partial_tp(t)\in H$
distributionally, so that
$T(t)\in H(\operatorname{div};\Omega)$.
The Green identity
\cite[Ch.~I, Sec.~2.1, Thm.~2.2 and Cor.~2.1]{GiraultRaviart1979}
defines its weak boundary traction
$\mathfrak t_{T(t)}\in\U$ by
\[
  \dualprod{\mathfrak t_{T(t)}}{\gamma_Pw}_{\U,G}
 =\scprod{\partial_tp(t)}{w}_H
 +\scprod{T(t)}{\nabla w}_{\Lp{2}(\Omega;\C^{d\times d})}
 \qquad(w\in W).
\]
By \eqref{kv:eq:weakpressure},
$\mathfrak t_{T(t)}=u(t)$ for the complete stress,
including pressure.

It remains to recover the displacement.
Since $v(t)\in V$, $v\in\conC(I;H)$,
and $Av=r\in\conC(I;\Hsym)$,
Korn's inequality gives $v\in\conC(I;V)$.
The descriptor equations also imply
$\partial_ts=Av=\varepsilon(v)$.
Hence, if $s(t_0)=\varepsilon(q_0)$
for some $q_0\in W$, then
$q(t):=q_0+\int_{t_0}^t v(\tau)\,\dd\tau$
satisfies $\partial_tq=v$ and
$s=\varepsilon(q)$ for all $t\in I$.

Conversely, suppose
\eqref{kv:eq:pde}--\eqref{kv:eq:boundary}
hold pointwise in time with
$q\in\conC^1(I;W)$,
$p\in\conC^1(I;H)$,
$u\in\conC(I;\U)$, and
$\pi(t)\in\Lp{2}(\Omega)$ for every $t\in I$.
Then $v=\partial_tq=\rho^{-1}p\in V$,
and the Green identity yields the domain condition for $C\adjun$
with elastic stress $2\mu\varepsilon(q)$, viscous stress
$2\eta Av$, and boundary traction $u$.
Thus a classical descriptor trajectory is obtained with
$s=\varepsilon(q)$, $r=Av$,
$y=\gamma_Pv$, and
$\lambda=(I-P_\sigma)\partial_tp$.
The multiplier $\lambda\in H_{\nabla,P}$
has a unique $\Hk{1}$-potential vanishing
on $\Gamma_P$, which need not coincide
with the physical pressure
$\pi\in\Lp{2}(\Omega)$.
The additive constant of the pressure is fixed by
the prescribed traction. By Korn's inequality,
$\varepsilon:W\to\Hsym$ is injective with closed range.
The constructions above are inverse to each other,
establishing a bijection between classical solutions
of \eqref{kv:eq:pde}--\eqref{kv:eq:boundary}
with the stated regularity and classical descriptor
trajectories satisfying
$s(t_0)\in\varepsilon(W)$.

The power balance \eqref{eq:classical-power}
takes the form \[
   \ddt\mathcal H(s(t),p(t))
  =-\int_\Omega2\eta|\varepsilon(v(t))|^2\dx[x]
  +\Re\dualprod{u(t)}{y(t)}_{\U,\U\adjun}.
\]

\backmatter
\section*{Statements and Declarations}

\bmhead{Funding}
This work was funded by the Deutsche Forschungsgemeinschaft
(DFG, German Research Foundation) -- Project-ID 531152215 -- CRC 1701.

\bmhead{Competing interests}
The authors declare that they have no competing interests.

\bmhead{Ethics approval and consent to participate}
Not applicable.

\bmhead{Consent for publication}
Not applicable.

\bmhead{Data availability}
No datasets were generated or analysed during this study.

\bmhead{Materials availability}
Not applicable.

\bmhead{Code availability}
Not applicable.

\bmhead{Authors' contributions}
Both authors contributed to the mathematical analysis and the
preparation of the manuscript.

\bibliography{Constraints}
\end{document}